\documentclass[12pt]{amsart} 

\usepackage{amssymb} 
\usepackage{amsmath} 
\usepackage{fullpage} 
\usepackage{url}

\newtheorem{theo}{Theorem}

\newtheorem{lemma}{Lemma}

\theoremstyle{definition}
\newtheorem{defi}{Definition}
\theoremstyle{remark}
\newtheorem{rem}{Remark}

\newtheorem{exer}{Exercise}
\def\Er{{\mathbb E}}

\def\Pr{{\mathbb P}}
\def\Qr{{\mathbb Q}}
\def\Rr{{\mathbb R}}

\def\Ac{{\mathcal{A}}}
\def\Bc{{\mathcal{B}}}

\def\Fc{{\mathcal{F}}}
\def\Gc{{\mathcal{G}}}

\def\Nc{{\mathcal{N}}}

\def\Sc{{\mathcal{S}}}

\def\one{{\rm \bf 1}}

\def\essinf{\operatorname{ess.\!inf}}

\def\({\left(}     
\def\){\right)}    
\def\[{\left[}     
\def\]{\right]}

\def\as{{\frenchspacing a.s.}~}
\begin{document}
\title{Commonotonicity and time consistency for Lebesgue continuous monetary utility functions } 
\author{Freddy Delbaen}
\address{Departement f\"ur Mathematik, ETH Z\"urich, R\"{a}mistrasse
101, 8092 Z\"{u}rich, Switzerland}
\address{Institut f\"ur Mathematik,
Universit\"at Z\"urich, Winterthurerstrasse 190,
8057 Z\"urich, Switzerland}
\date{First version December 2018, this version \today}

\begin{abstract}
It is proved that commonotonicity and time consistence for monetary utility functions do not go together. I also gives additional results on atomless and conditionally atomless probability spaces.\end{abstract}

\maketitle

\section{Notation}

The purpose of this paper is to investigate the relation between commonotonicity and time consistency of monetary utility functions.\footnote{ This research was done while the author was visiting Tokyo Metropolitan University in October and November 2018. I thank the staff of TMU for the many fruitful discussions and in particular I thank Prof. Adachi for the many critical remarks.}  It will turn out that it is sufficient to have a two period model.  In this setting we will work with a probability space equipped with three sigma algebras $(\Omega,\Fc_0\subset \Fc_1\subset \Fc_2,\Pr)$.  The sigma algebra $\Fc_0$ will be supposed to be trivial $\Fc_0=\{\emptyset,\Omega\}$ whereas the sigma algebra $\Fc_2$ will be supposed to express innovations with respect to $\Fc_1$.  Since we do not put topological properties on the set $\Omega$ we will make precise definitions later that do not use conditional probability kernels.  But essentially we could say that we suppose that conditionally on $\Fc_1$ the probability $\Pr$ is atomless on $\Fc_2$. We will show that such an hypothesis implies that there is an atomless sigma algebra $\Bc\subset \Fc_2$ that is independent of $\Fc_1$.   The space $L^\infty(\Fc_i)$ will denote the space of bounded $\Fc_i$ measurable random variables, modulo almost sure equality \as. We will also suppose that there is a time consistent utility function $u_2\colon L^\infty(\Fc_2)\rightarrow \Rr$.  As shown in \cite{FDbook} this means that we also have utility functions $u_{1,2}\colon L^\infty(\Fc_2)\rightarrow L^\infty(\Fc_1)$ and $u_{0,1}\colon L^\infty(\Fc_1)\rightarrow L^\infty(\Fc_0)=\Rr$ such that $u_2=u_{0,1}\circ u_{1,2}$.  In particular $u_{0,1}$ is simply the restriction of $u_2$ to $L^\infty(\Fc_1)$.  Our utility functions are monetary and concave which is expressed in the following list of properties
\begin{enumerate}
\item  For $i<j$ we have $u_{i,j}\colon L^\infty(\Fc_j)\rightarrow L^\infty(\Fc_i)$, if $\xi\ge 0$ then also $u_{i,j}(\xi)\ge 0$ and $u_{i,j}(0)=0$.
\item For $\xi,\eta\in L^\infty(\Fc_j)$, $0\le\lambda\le1$ and $\Fc_i$ measurable we have
$$
u_{i,j}(\lambda \xi +(1-\lambda)\eta)\ge \lambda u_{i,j}(\xi)+(1-\lambda u_{i,j}(\eta).
$$
Since commonotonicity implies positive homogeneity we will use a stronger property and suppose coherence:
\item For $\xi\in L^\infty(\Fc_j)$, $0\le\lambda$ and $\Fc_i$ measurable we have
$$
u_{i,j}(\lambda \xi)=\lambda u_{i,j}(\xi).
$$
\item For $\xi\in L^\infty(\Fc_j)$ and $a\in L^\infty(\Fc_i)$ we have
$$
u_{i,j}(\xi +a)=u_{i,j}(\xi) + a.
$$
\item We will need Lebesgue continuity which means:  if $\xi_n\in L^\infty(\Fc_j)$ is a uniformly bounded sequence such that $\xi_n\rightarrow \eta$ in probability then also $u_{i,j}(\xi_n)$ tends to $u_{i,j}(\eta)$ in probability.
\end{enumerate}
The utility functions we need are coherent and hence we can use the dual representation.  There is a set of probability measures, $\Sc$, absolutely continuous with respect to $\Pr$ such that
$$
u(\xi)=\inf\Er_{\Qr\in \Sc}[\xi].
$$
The set is seen as a subset of $L^1$ and is supposed to be convex and closed.  The Lebesgue continuity is equivalent to the weak compactness of $\Sc$.  We will suppose that our utility functions are relevant,  i.e. for each $A$ with $\Pr[A]>0$ we have $u(-\one_A)<0$. see \cite{FDbook}.  By the Halmos-Savage theorem this means that $\Sc$ contains an equivalent probability measure.  We need this property in order to avoid some problems with negligible sets appearing in the definition and comparison of conditional expectations.

We say that two random variables $\xi,\eta$ are commonotone if there are two nondecreasing functions $f,g\colon \Rr\rightarrow \Rr$ and a random variable $\zeta$ such that $\xi=f(\zeta), \eta=g(\zeta)$.  {\it Commonotonicity can be seen as the opposite of diversification}.  If $\zeta$  increases then both $\xi$ and $\eta$ increase (or better do not decrease). By the way in case $\xi$ and $\eta$ are commonotone then one can choose $\zeta=\xi+\eta$, see \cite{FDbook}.  It can be shown (an exercise) that in this case one can choose representatives (still denoted $\xi,\eta$ such that $(\xi(\omega)-\xi(\omega'))(\eta(\omega)-\eta(\omega'))\ge 0$ for all $\omega,\omega'$.  We say that a set $E\subset \Rr^2$ is commonotone  if $(x,y),(x',y')\in E$ implies $(x-x')(y-y')\ge 0$.  Using this, random variables $\xi,\eta$ are commonotone if and only if, the support of the image measure of $(\xi,\eta)$ is a commonotone set. We will need the following obvious result
\begin{lemma}  If $(\alpha,\beta)\in \Rr^2$ and the couple $(\xi,\eta)\colon\Omega\rightarrow\Rr^2$ takes values in the set
$$
\left\{ (x,\beta)\mid x\le \alpha   \right\} \cup \left\{ (\alpha,y)\mid y\ge \beta   \right\},
$$
then $\xi$ and $\eta$ are commonotone.
\end{lemma}
The concept of time consistency (and inconsistency) was already investigated by Koopmans, \cite{Koop1}. The role of commonotonicity found its way in insurance and is present in several papers.  The use of Choquet integration as premium principle was emphasized by Denneberg, \cite{Denn}. Denneberg was inspired by the pioneering work of Yaari, \cite{Yaa}. Schmeidler proved the relation between commonotone principles, convex games and Choquet integration, \cite{Schm}.  Modern uses can be found in for instance \cite{WYP} and \cite{Wang}.  The concept of risk measures (up to sign changes monetary utility functions) was introduced in \cite{ADEH1} and \cite{ADEH2}.  For more references and different proofs of these results I refer to \cite{FDbook}.

\section{Atomless Extension}

\begin{defi}  We say that $\Fc_2$ is atomless conditionally to $\Fc_1$ if the following holds.  If $A\in \Fc_2$ then there exists a set $B\subset A$, $B\in\Fc_2$, such that $0< \Er[\one_B\mid\Fc_1]<\Er[\one_A\mid\Fc_1]$ on the set $\{\Er[\one_A\mid\Fc_1]>0\}$.
\end{defi}
In  case the conditional expectation could be calculated with a -- under extra topological conditions -- regular probability kernel, say $K(\omega, A)$, then the above definition is a measure theoretic way of saying that the probability measure $K(\omega, .)$ is atomless for almost every $\omega\in\Omega$.
\begin{rem} Using an exhaustion argument, one can show that the definition is equivalent with the existence of $B\subset A$ such that $\Pr\[0< \Er[\one_B\mid\Fc_1]<\Er[\one_A\mid\Fc_1]\]>0$.
\end{rem}
The main result of this section is the following
\begin{theo} $\Fc_2$ is atomless conditionally to $\Fc_1$ if and only if there exists an atomless sigma algebra $\Bc\subset \Fc_2$ that is independent of $\Fc_1$.
\end{theo}
The ``if" part is easy but requires some continuity argument.  Because $\Bc$ is atomless, there is a $\Bc$-measurable, $[0,1]$ uniformly distributed random variable $U$. The sets $B_t=\{U\le t\}, 0\le t \le 1$ form an increasing family of sets with $\Pr[B_t]=t$. Let $A\in \Fc_2$  and let $F=\{ 0 < \Er[\one_A\mid \Fc_1]\}$.  We may suppose that $\Pr[F]>0$ since otherwise there is nothing to prove. We will show that there is $t\in]0,1[$ with $\Pr\[ 0 < \Er[\one_{A\cap B_t} \mid \Fc_1]< \Er[\one_A\mid \Fc_1]    \] > 0$. Obviously for $0\le s\le t \le 1$ we have, by independence of $\Bc$ and $\Fc_1$:
$$
\Vert   \Er[\one_{A\cap B_t} \mid \Fc_1] -  \Er[\one_{A\cap B_s} \mid \Fc_1]\Vert_\infty \le \Vert   \Er[\one_{B_t\setminus B_s} \mid \Fc_1]\Vert_\infty = t-s.
$$
It follows that there is a set of measure $1$, say $\Omega'$,  such that for all $s\le t$, rational, $$| \Er[\one_{A\cap B_t} \mid \Fc_1] -  \Er[\one_{A\cap B_s} \mid \Fc_1] | \le t-s$$ on $\Omega'$.  On the set $\Omega'$ we can extend these functions $$\{q\in [0,1] \mid q \text{ rational }\}\rightarrow \Er[\one_{A\cap B_q}\mid \Fc_1](\omega)$$ to a continuous function on $[0,1]$.  The resulting continuous extension then represents $\(\Er
[\one_{A\cap B_t} \mid \Fc_1]\)_{t}$. For $t=0$ we have zero and for $t=1$ we find $\Er[\one_A\mid \Fc_1]$. Because for $\omega\in\Omega'$, the trajectories are continuous we must by a simple application of Fubini's theorem, have that the real valued function
$$
t\rightarrow \Pr\[ 0 < \Er[\one_{A\cap B_t} \mid \Fc_1]< \Er[\one_A\mid \Fc_1]    \] 
$$ 
becomes strictly positive for some $t$.  According to remark 1, $\Fc_2$ is conditionally atomless with respect to $\Fc_1$. With some extra work -- as will be done later -- one can even show that there is $G\subset A$ such that $\Er[\one_G\mid\Fc_1]= (1/2)\Er[\one_A\mid\Fc_1]$.

The proof of the ``only if" part is broken down in several steps.  We will without further mentioning, always suppose that $\Fc_2$ is atomless conditionally to $\Fc_1$.
\begin{lemma} Suppose $A\in\Fc_1$ and $C\subset A$ is such that $\Er[\one_C\mid \Fc_1]>0$ on $A$.  Then we can construct a decreasing sequence of sets $(B_n)_{n\ge 0}$, $B_n\subset C$, such that  $0<\Er[\one_{B_n}\mid\Fc_1]\le 2^{-n}$ on $A$.
\end{lemma}
{\bf Proof}  The statement is obviously true for $n=0$ since we can take $B_0=C$.  We now proceed by induction and suppose the statement holds for $n$. So the set $B_n\subset A$ satisfies $0<\Er[\one_{B_n}\mid\Fc_1]\le 2^{-n}$ on $A$.  Clearly $\{\Er[\one_{B_n}\mid\Fc_1]>0\}=A$.  By assumption there is a set $D\subset B_n$ such that on $A$ we have
$$
0<\Er[\one_D\mid\Fc_1]<\Er[\one_{B_n}\mid\Fc_1].
$$
We now take
$$
B_{n+1}=\(D\cap \left\{\Er[\one_D\mid\Fc_1]\le \frac{1}{2}\Er[\one_{B_n}\mid\Fc_1]\right\}\)\cup
\((B_n\setminus D)\cap\left\{\Er[\one_D\mid\Fc_1]> \frac{1}{2}\Er[\one_{B_n}\mid\Fc_1]\right\}\).
$$
The set $B_{n+1}$ satisfies the requirements.
\begin{lemma} Let $C\in\Fc_2$ and let $h\colon\Omega\rightarrow [0.1]$ be $\Fc_1$ measurable.  Then there is a set $B\subset C$ such that $\Er[\one_B\mid\Fc_1]=h\,\Er[\one_C\mid\Fc_1]$.
\end{lemma}
{\bf Proof}  Let $\Bc$ be the class
$$
\Bc=\left\{ B\subset C\mid \Er[\one_B\mid\Fc_1]\le  h\,\Er[\one_C\mid\Fc_1]  \right\}.
$$
Let $B_\alpha$ be a totally ordered system where $\alpha$ runs through some set $I$.  Let $m=\sup_\alpha \Pr[B_\alpha]$.  There is an increasing sequence of sets $B_{\alpha_n}$ such that $\Pr[B_{\alpha_n}]\uparrow m$.  Let $B_\infty=\cup_n B_{\alpha_n}$.  Obviously, by monotone convergence,
$B_\infty\in \Bc$.  Because $\Pr[B_\infty]=m$ is it easily seen that we must have $B_\alpha\subset B_\infty$ for every $\alpha\in I$.  Because every totally ordered system in $\Bc$ has a majorant, the class $\Bc$ must have a maximal element.  (This is the usual Zorn's lemma).  Let $B\in\Bc$ be a maximal element.  We already have
$$
\Er[\one_B\mid\Fc_1]\le  h\,\Er[\one_C\mid\Fc_1].
$$
Suppose the $\Pr[\Er[\one_B\mid\Fc_1] <  h\,\Er[\one_C\mid\Fc_1] ]>0$.  Then there is $n$ such that $\Pr[\Er[\one_B\mid\Fc_1] <  h\,\Er[\one_C\mid\Fc_1] - 2^{-n}]>0$. Also $\Pr[C\setminus B]>0$.  The previously lemma then allows to find a set $D\subset C\setminus B$ with $0<\Er[\one_D\mid\Fc_1]\le 2^{-n}$ on the set $\{\Er[\one_B\mid\Fc_1] <  h\,\Er[\one_C\mid\Fc_1] - 2^{-n}\}$ and zero elsewhere.  The set $B\cup D$ is still in $\Bc$ and is strictly bigger than the maximal elemet $B$, a contradiction.  This shows that
$$
\Er[\one_B\mid\Fc_1] = h\,\Er[\one_C\mid\Fc_1],
$$
as requested.
\begin{theo} There is an increassing family of sets $(B_t)_{t\in [0,1]}$ such that $\Er[\one_{B_t}\mid\Fc_1]=t$.  The sigma algebra $\Bc$, generated by the family $(B_t)_t$ is independent of $\Fc_1$.  The system $(B_t)_t$ can also be described as $B_t=\{U\le t\}$ where $U$ is a random variable that is uniformly distributed on $[0,1]$,  $U$ and $\Fc_1$ are independent.
\end{theo}
{\bf Proof} The proof is a repeated use of the previous lemma where we take $h=1/2$. We start with $B_0=\emptyset, B_1=\Omega$.  Suppose that for the diadic numbers $k 2^{-n}, k=0,\ldots 2^n$ the sets are already defined.  Then we consider the set $B_{(k+1)2^{-n}}\setminus B_{k2^{-n}}$ and apply the previous lemma with $h=1/2$.  We get a set $D\subset B_{(k+1)2^{-n}}\setminus B_{k2^{-n}}$ with $\Er[\one_D\mid\Fc_1]=2^{-(n+1)}$.  We then define $B_{(2k+1)2^{-(n+1)}}=B_{k2^{-n}}\cup D$.  For non-diadic numbers $t$ we find a sequence of diadic numbers $d_n$ such that $d_n\uparrow t$.  Then we define $B_t=\cup_n B_{d_n}$.  This completes the construction.  Since the system $(B_t)_t$ is trivially stable for intersection, the relation $\Er[\one_{B_t}\mid\Fc_1]=t$ shows that the sigma algebra $\Bc$ generated by $(B_t)_t$ is independent of $\Fc_1$.  The construction of $U$ is standard.  At level $n$ we put $U_n=\sum_{k=1,\ldots 2^n} k2^{-n}\one_{B_{k2^{-n} }\setminus B_{(k-1)2^{-n}}}$.  $U_n$ then decreases to a random variable $U$ that satisfies the needed properties.
\begin{theo}
For every  $\Fc_1$ measurable function $h\colon \Omega\rightarrow [0,1]$, there is a set $B_h\in\Fc_2$ such that $\Er[\one_{B_h}\mid\Fc_1]=h$.
\end{theo}
{\bf Proof} The idea is to use the set $B_t$ on the set $\{h=t\}$, i.e. $B=\cup_t(\{h=t\}\cap B_t)$.  However, because the set of real numbers is uncountable, this definition is not good enough to obtain a set in $\Fc_2$. So we need a trick.  Let $\phi$ be the mapping
$$
\phi\colon (\Omega,\Fc_2)\rightarrow (\Omega,\Fc_1)\times(\Omega,\Bc), \phi(\omega)=(\omega,\omega).
$$
This mapping is obviously measurable and the image measure is -- because of independence -- the product measure.  We also define $h_1(\omega,\omega')=h(\omega)$ and $U_2(\omega,\omega')=U(\omega')$.   For  $A\in \Fc_1$ we denote $A_1=A\times\Omega$.  Now we put $B_h=\{U\le h\}=\phi^{-1}\{U_2\le h_1\}$.  We now verify that $\Er[\one_{B_h}\mid\Fc_1]=h$.  To do this we calculate for a set $A\in\Fc_1$ the probability $\Pr[B_h\cap A]$.
\begin{eqnarray*}
\Pr[B_h\cap A] &=& \Pr\times\Pr[(U_2\le h_1)\cap A_1] \\
&=& \int \Pr[d\omega']\int \Pr[d\omega] \one_{\{U_2\le h_1\} }(\omega,\omega')\one_{A_1}(\omega,\omega')\\
&=& \int\Pr[d\omega']\Pr[\{h\ge U(\omega'\}\cap A]\\
&=& \int_0^1 dt\, \Pr[\{h\ge t\}\cap A]\\
&=& \Er[h\one_A],
\end{eqnarray*}
showing $\Er[\one_{B_h}\mid\Fc_1]=h$.
\begin{rem}  The previous theorem is not actually needed.  We will need the stronger version where the conditional expectation is replaced by the utility function $u_{1,2}$.  To prove this stronger version we will use a slightly different approach.  However in the case where we are only interested in conditional expectations the above proof might be of some didactical interest.
\end{rem}
\begin{rem} After the  first version was made available, I got the remark that the paper \cite{SSWW} of Shen, J., Shen, Y.,  Wang, B., and  Wang, R.  contains similar concepts and results.\footnote{I thank Ruodu Wang for pointing out these relations and for the subsequent discussions we had on the topic} In their notation they work with a measurable space $(\Omega,\Ac)$ on which they have a finite number of probability measures $\Qr_1,\ldots,\Qr_n$.\footnote{Their paper also considers an infinite number of measures but to clarify the relation between their paper and my approach, I only consider a finite number of measures} They then say
\begin{defi} The set $(\Qr_1,\ldots,\Qr_n)$ is conditionally atomless if there exists a dominating measure $\Qr$ (i..e $\Qr_k\ll \Qr$ for each $k\le n$) as well as a continuously distributed random variable $X$ (for the measure $\Qr$) such that the vector of Radon-Nikodym derivatives $\(\frac{d\Qr_k}{d\Qr}\)_{k}$ is independent of $X$.
\end{defi}
They then prove the following
\begin{theo} Are equivalent
\begin{enumerate}
\item $(\Qr_1,\ldots,\Qr_n)$ is conditionally atomless
\item in the definition we can take $\Qr=\frac{1}{n}(\Qr_1+\ldots+\Qr_n)$
\item $X$ can be taken as uniformly distributed over $[0,1]$.
\end{enumerate}
\end{theo}
There are several differences with my approach.  There is the technical differenc that they suppose the existence of a continuously distributed random variable $X$.  In doing so they avoid the technical points between the more conceptual definition using conditional expectations and the construction of a suitable sigma-algebra with a uniformly distributed random variable.   A further difference is that they use a dominating measure that later can be taken as the mean of $(\Qr_1,\ldots,\Qr_n)$.  Of course their result together with the results here show that the definition that $(\Qr_1,\ldots,\Qr_n)$ is conditionally atomless is equivalent to the statement that for the measure $\Qr_0=\frac{1}{n}(\Qr_1+\ldots+\Qr_n)$, the sigma algebra $\Ac$ is conditionally atomless with respect to the sigma-algebra generated by the Radon-Nikodym derivatives $\(\frac{d\Qr_k}{d\Qr_0}\)_k$.  In \cite{SSWW} it is also shown that one can take any strictly positive convex combination of the measures $(\Qr_1,\ldots,\Qr_n)$. Below we will show that this sigma-algebra in some sense has a minimal property, a result that clarifies the relation between the two approaches.  Before doing so, let us recall two easy exercises from introductory probability theory.
\begin{exer}  For a given probability space $(\Omega,\Ac,\Qr)$ let us denote $\Nc=\{N\in\Ac\mid\Qr[N]=0\}$.  Suppose that a sub sigma-algebra $\Fc_\subset\Ac$ is given and that $\Gc$ is another sub sigma-algebra which is included in the sigma-algebra generated by $\Fc$ and $\Nc$. Then for each $\xi\in L^1(\Omega,\Ac,\Qr)$
$$
\Er_\Qr[\xi\mid \Fc]=\Er_\Qr[\xi\mid \Gc] \quad\text{\as}
$$
\end{exer}
\begin{exer} With the notation in the previous exercise let $F\colon\Omega\rightarrow \Rr^n$ and $F'\colon\Omega\rightarrow \Rr^n$ be two vectors that are equal \as.  Let $\Fc$ be generated by $F$ and $\Gc$ be generated by $F'$.  Then $\Fc$ and $\Gc$ are equal up to sets in $\Nc$.  More precisely $\Gc$ is included in the sigma-algebra generated by $\Fc$ and $\Nc$ (and of course conversely), i.e. $\sigma(\Fc,\Nc)=\sigma(\Gc,\Nc)$.
\end{exer}
\begin{theo} Let $\Qr_1,\dots,\Qr_n$ be probability measures on a measurable space $(\Omega,\Ac)$.  Let $\Qr_0$ denote  a convex combination of these measures $\Qr_0=\sum_k \lambda_k \Qr_k$ where each $\lambda_k > 0$. Let $f_k$ denote an $\Ac$ measurable version $\frac{d\Qr_k}{\Qr_0}$.  Let $\Qr$ be another dominating measure with $g_k$ an $\Ac$ measurable version of $\frac{d\Qr_k}{d\Qr}$.  Let $\Nc=\{N\in \Ac\mid \Qr_0[N]=0\}$.  Let $\Fc$ be generated by $f_k,k=1\dots n$ and $\Gc$ be generated by $g_k,k=1\ldots n$.  Then $\Fc\subset \sigma(\Gc,\Nc)$
\end{theo}
{\bf Proof }  Clearly $\Qr_0\ll \Qr$ so let $h=\frac{d\Qr_0}{d\Qr}$. It is now immediate that $g_k= f_k h$ $\Qr$ \as. To see this, observe that the values of $f_k$ on $\{h=0\}$ do not matter. The functions $g_k$ and $h$ are $\Gc$ measurable since $h$ can be taken as $h=\sum_k \lambda_k g_k$. Then we define $f_k'= \frac{g_k}{h}$ on $\{h>0\}$ and $f_k'=0$ on $\{h=0\}$. This choice shows that the $f_k'$ are $\Gc$ measurable. It is immediate that $f_k=f_k'$ $\Qr_0$ \as.   The result now follows from the exercises.

From the theorem it follows that the sigma-algebra augmented with the class $\Nc$ is the same for all strictly positive convex combinations.  The exercise shows that in the definition of conditionally atomless with respect to $\Fc$, we can also add the null sets $\Nc$ to $\Fc$.  To check that $\Ac$ is conditionally atomless with respect to a sigma-algebra $\Fc$ it is clear that the smaller $\Fc$, the easier it is to satisfy the condition. In my opinion the above clarifies the relation between this paper and \cite{SSWW}.
\end{rem}
\section{A continuity Result}
For each $h\colon \Omega\rightarrow [0,1]$ that is $\Fc_1$ measurable we put $\phi(h)=u_{1,2}(\one_{\{U\le h\}})$. Clearly $\phi$ takes values in the space $L^\infty(\Fc_1)$.  We have the following continuity result.
\begin{theo} Suppose that $u_{1,2}$ has the Fatou property and the Lebesgue property. Suppose that conditional;ly to $\Fc_1$, the sigma algebra $\Fc_2$  is atomless. If $h_n\downarrow h$ or $h_n\uparrow h$ we have $\phi(h_n)\rightarrow \phi(h)$.
\end{theo}
{\bf Proof }  If $h_n\downarrow h$ then $\one_{\{U\le h_n\}}\downarrow \one_{\{U\le h\}}$ and the Fatou property gives the desired result. For the upward convergence we must be more careful.  Because $U$ has a continuous distribution function and since it is independent of $\Fc_1$, we conclude that $\Pr[U=h]=0$ and hence 
$\one_{\{U\le h_n\}}\uparrow \one_{\{U\le h\}}$ \as .  The Lebesgue property then allows to conclude.

\begin{theo} If $h\colon\Omega\rightarrow [0,1]$ is $\Fc_1$ measurable, there is an $\Fc_1$ measurable function $g\colon\Omega\rightarrow [0,1]$ such that  set $B_g=\{U\le g\}$ satisfies $u_{1,2}(\one_{B_g})=h$.
\end{theo}
{\bf Proof} The statement can be rewritten as $\phi(g)=h$. Let us introduce the class
$$
\Gc=\{ g\mid g\text{ is } \Fc_1\text{ measurable and }u_{1,2}(\one_{B_g})=\phi(g)\ge h  \}.
$$
$\Gc$ is nonempty since $1\in\Gc$.  Furthermore $\Gc$ is stable for taking the minimum.  Indeed, let $g_1,g_2\in \Gc$ and put $g=g_1\one_A+g_2\one_{A^c}$ where $A=\{g_1<g_2\}$.  Since $u_{1,2}(g)=u_{1,2}(g_1)\one_A + \one_{A^c}u_{1,2}(g_2)\ge h$ we have that $g\in\Gc$.  Let now $g_n\downarrow g$ where $g_n\in\Gc$ and $\Er[g_n]\downarrow \inf\{\Er[g']\mid g'\in\Gc\}$.  The continuity for decreasing sequences then shows that $g\in\Gc$.  The continuity for increasing sequences (the Lebesgue property) will show that actually $u_{1,2}(\one_{B_g})=h$.  Suppose on the contrary that the set $\{u_{1,2}(\one_{B_g})>h\}$ has non zero measure.  Take now a sequence $g_n\uparrow g$ such that $g_n<g$ almost everywhere. By the previous theorem $u_{1,2}(\one_{B_{g_n}})\uparrow u_{1,2}(\one_{B_g})$.  Hence, there must exist $n$ such that $A_n=\{u_{1,2}(\one_{B_{g_n}})>h\}$ has non zero measure.  Put now $g'=g_n\one_{A_n}+g\one_{A_n^c}$.  We have $\Er[g']<\Er[g]$ but also $g'\in \Gc$ a contradiction to the minimality of $g$.
\begin{rem}  Although ``intuitively clear", the continuity of the process $u_{1,2}(\one_{B_t})$ is not an easy result.  First of all, we are working with random variables identified under the equivalence \as. That means that we must first select or construct measurable functions instead of classes of measurable functions.  Then we must show that with respect to $t$ these outcomes are continuous.  The general theory of stochastic processes gives us the necessary tools to achieve this goal.  We do not really need these finer results so the remark can be skipped if you do not belong to the amateurs of the general theory of stochastic processes, see \cite{DM} for the necessary details.  First we will construct a process $\alpha(t,\omega)$. For each rational point $q\in [0,1]$ we select an $\Fc_1$ measurable function $\alpha'(q)$ that represents $u_{1,2}(\one_{B_q})$.  Because of monotonicity we can -- if needed -- change these selections on a set of zero measure, to make sure that \as the mapping $\Qr\cap [0,1]\rightarrow \Rr; q\rightarrow \alpha'(q)$ is increasing.   For each $t\in [0,1]$ we now define $\alpha(t)=\inf_{q\text{ rational },q\ge t}\alpha'(q)$.  The functions $\alpha(t)$ are of course $\Fc_1$ measurable and represent $u_{1,2}(\one_{B_t})$ by the Fatou property.   We may also suppose that  $\alpha(0)=0,\alpha(1)=1$ \as. It is clear that $\alpha$ is \as increasing in $t$ and is right continuous.  This means there is a set (independent of $t$) such that on this set $t\rightarrow\alpha(t,\omega)$ is right continuous and increasing.

The function $\alpha$ also satisfies $\alpha(h)=u_{1,2}(\one_{\{U\le h\}})=\phi(h)$ for each $\Fc_1$ measurable function $h\colon \Omega\rightarrow [0,1]$.  \footnote{ To avoid misunderstandings the random variable $\alpha(h)$ is defined as $\alpha(h)(\omega)=\alpha(h(\omega),\omega)$. Such a practice is common in stochastic process theory.}  The statement is easy to verify for elementary functions $h$ and the general statement trivially follows by approximating $h$ from {\it above} by elementary functions.  Let us give the details.  For en elementary function $h=\sum_{k=1}^K t_k\one_{A_k}$ (the sets $A_k$ are disjoint and in $\Fc_1$), we have
\begin{eqnarray*}
\alpha(h) &=& \sum_k \alpha(t_k)\one_{A_k}\\
&=&\sum_k u_{1,2}(\one_{B_{t_k}})\one_{A_k}\\
&=&\sum_k u_{1,2}(\one_{B_{t_k}}\one_{A_k})\one_{A_k}\\
&=&\sum_k u_{1,2} \( \one_{B_{t_k}\cap A_k}\)\one_{A_k}\\
&=&\sum_k u_{1,2}\(\(\sum_l \one_{B_{t_l}\cap A_l}\)\one_{A_k} \)\one_{A_k}\\
&=&\sum_k u_{1,2}(\one_{\{U\le h\}}\one_{A_k})\one_{A_k}\\
&=& u_{1,2}(\one_{\{U\le h\}})=\phi(h).
\end{eqnarray*}
As indicated above the Fatou property then completes the proof using right continuity.  Indeed, let $h\colon \Omega\rightarrow [0,1]$ be $\Fc_1$ measurable and let $h_n \downarrow h$ be a sequence of elementary functions, that are $\Fc_1$ measurable.  Since $\one_{\{U\le h_n\}}
\downarrow \one_{\{U\le h\}}$, the Fatou property and the right continuity of $\alpha(t)$ give us $\phi(h)=u_{1,2}(\one_{\{U\le h\}})$. 

The proof of the left continuity can be done using the general theory of stochastic processes.  The basic ingredient is the theorem in the beginning of this section. We leave the details as an exercise.  To conclude this remark, we can say that the process $\alpha_t=\phi(t)$ has increasing and continuous trajectories.
\end{rem}
\section{The main Result}
This section is devoted to the proof of
\begin{theo}  Suppose that $u_{1,2}$ is Lebesgue continuous and suppose that $\Fc_2$ is atomless conditionally to $\Fc_1$. For any two bounded $\Fc_1$ measurable random variables, $f,g$, we can find two commonotone $\Fc_2$ random variables $\xi,\eta$ such that $f=u_{1,2}(\xi), g=u_{1,2}(\eta)$ and $u_{1,2}(\xi+\eta)=f+g$.  The random variables also satisfy $\Vert\xi\Vert_\infty,\Vert\eta\Vert_\infty\le 3 \max(\Vert f\Vert_\infty,\Vert g\Vert_\infty)$.
\end{theo}
Let $m=\max(\Vert f\Vert_\infty,\Vert g\Vert_\infty)$.  We will construct the two random variables so that they take values in the commonotone set
$$
V=\left\{ (x,-m)\mid x\le m   \right\} \cup \left\{ (m,y)\mid y\ge -m   \right\}.
$$
As observed before, these two random variables are then commonotone.  We remark that the couple $(f,g)$ takes its values in the set
$$
W=\left\{(x,y)\mid x\le m, y\ge -m \right\}.
$$
For $(x,y)\in W$ we will define three functions $X\colon W\rightarrow V, Y\colon W\rightarrow V, \lambda\colon W\rightarrow [0,1]$.  For $(x,y)\in W, (x,y)\neq (m,-m)$ we look at the two points of intersection of the line $\{(x+t,y+t)\mid t\in\Rr\}$ and $V$.  The intersection with the halfline $\left\{ (x,-m)\mid x\le m   \right\}$ is $X$, the intersection with the halfline $\left\{ (m,y)\mid y\ge -m   \right\}$ is $Y$.  The function $\lambda$ is defined so that $\lambda Y + (1-\lambda) X=(x,y)$.  At the point $(m,-m)$ we get $X=Y=(x,y)$ and we put $\lambda=0$.  The functions $X$ and $Y$ are continuous on $W$ whereas $\lambda$ is only continuous on $W\setminus\{(m,-m)\}$. We can write $(x,y)=X+\lambda (Y-X)$.  We remark that both coordinates of $Y-X$ are positive.  For the given functions $f,g$ we now consider the composite functions $X(f,g),Y(f,g),\lambda(f,g)$.  All these functions are $\Fc_1$ measurable.  For $\lambda(f,g)$ we take a set $B\in\Fc_2$ such that $u_{1,2}(\one_B)=\lambda(f,g)$.  The existence of $B$ follows from the results of the previous section.  The functions $\xi,\eta$ are now defined by the relation $(\xi,\eta)=X(f,g)+\one_B(Y(f,g)-X(f,g))=\one_{B^c}X(f,g)+\one_B Y(f,g)$.  Because of positive homogeneity we get $u_{1,2}(\xi)=f; u_{1,2}(\eta)=g$.  Since $(\xi,\eta)$ takes values in $V$, the two coordinates are commonotone. The equality $u_{1,2}(\xi+\eta)=f+g$ again follows  from positive homogeneity and $u_{1,2}(\one_B)=\lambda(f,g)$. Elementary algebra shows that $\Vert \xi\Vert_\infty, \Vert\eta\Vert_\infty\le 3m$.
\section{Commonotonicity and Time Consistency}
In his section we use the same hypothesis on the filtration $(\Fc_0,\Fc_1,\Fc_2)$.  In particular we suppose that $\Fc_2$ is atomless conditionally to $\Fc_1$. We start with a monetary coherent utility function $u_{0,2}\colon L^\infty(\Fc_2)\rightarrow\Rr$.  We  suppose  -- as in the rest of the paper -- that $u_{0,2}$ is relevant. As shown in \cite{FDbook},  there is a way to check whether this utility function can be extended to a time consistent utility function.  To do this we introduce the acceptability cones $\Ac_{0,2}=\{\xi\mid u_{0,2}(\xi)\ge 0\}$, $\Ac_{0,1}=\{\xi\in L^\infty(\Fc_1)\mid u_{0,2}(\xi)\ge 0\}$,$\Ac_{1,2}=\{\xi\in L^\infty(\Fc_2)\mid\text{ for all } A\in \Fc_1:\, u_{0,2}(\xi \one_A)\ge 0\}$.  The necessay and sufficient condition for the existence of a time consistent extension is $\Ac_{0,2}=\Ac_{0,1}+\Ac_{1,2}$.  If this is fulfilled we put $u_{1,2}(\xi)=\essinf \{\eta\in L^\infty(\Fc_1)\mid \xi-\eta\in\Ac_{1,2}\}$.  $u_{0,1}$ is simply the restriction of $u_{0,2}$ to $L^\infty(\Fc_1)$.
\begin{theo} Suppose
\begin{enumerate}
\item $u_{0,2}$ is coherent and relevant
\item $\Fc_2$ is atomless conditionally to $\Fc_1$
\item $u_{0,2}$ is time consistent
\item $u_{0,2}$ is commonotone, i.e. if $\xi,\eta\in L^\infty(\Fc_2)$ are commonotone, then $u_{0,2}(\xi+\eta)=u_{0,2}(\xi)+u_{0,2}(\eta)$
\item $u_{0,2}$ is Lebesgue continuous.
\end{enumerate}
Then there is a probability $\Qr\sim\Pr$ such that for all $f\in L^\infty(\Fc_1)$ we have $u_{0,1}(f)=\Er_\Qr[f]$.
\end{theo}
{\bf Proof }  According to the previous section for each $f,g\in L^\infty(\Fc_1)$ there are {\it commonotone } $\xi,\eta\in L^\infty(\Fc_2)$ with $u_{1,2}(\xi)=f,\, u_{1,2}(\eta)=g$ and $u_{1,2}(\xi+\eta)=f+g$. We then have $u_{0,1}(f)=u_{0,1}\(u_{1,2}(\xi)\)=u_{0,2}(\xi)$ and similarly for $g$.  The combination with commonotonicity then gives 
\begin{eqnarray*}
u_{0,1}(f+g)&=&u_{0,1}(u_{1,2}(\xi+\eta))\\
&=&u_{0,2}(\xi+\eta)\\
&=&u_{0,2}(\xi)+u_{0,2}(\eta)\\
&=&u_{0,1}(u_{1,2}(\xi))+u_{0,1}(u_{1,2}(\eta))\\
&=&u_{0,1}(f)+u_{0,1}(g)
\end{eqnarray*}
This shows that $u_{0,1}$ is additive (therefore linear) and hence is given by a finitely additive probability measure.  But Lebesgue continuity implies that this measure, say $\Qr$, should be sigma additive and absolutely continuous with respect to $\Pr$. Because $u_{0,2}$ and hence $u_{0,1}$ are relevant we must have $\Qr\sim\Pr$.
\begin{rem} For $\xi,\eta$ commonotone (and not just for the ones used in the proof of the theorem)  we can now  prove  that $u_{1,2}(\xi+\eta)=u_{1,2}(\xi)+u_{1,2}(\eta)$.  In fact this  holds for the equality $\Qr$ \as We already know that $u_{1,2}(\xi+\eta)\ge u_{1,2}(\xi)+u_{1,2}(\eta)$. If $\Qr[u_{1,2}(\xi+\eta)> u_{1,2}(\xi)+u_{1,2}(\eta)] >0$, then  we have 
\begin{eqnarray*}
u_{0,2}(\xi+\eta)&=&u_{0,1}(u_{1,2}(\xi+\eta))\\
&>& u_{0,1}(u_{1,2}(\xi)+u_{1,2}(\eta))\\
&\ge& u_{0,1}(u_{1,2}(\xi))+u_{0,1}(u_{1,2}(\eta))\\
&=& u_{0,2}(\xi)+u_{0,2}(\eta)
\end{eqnarray*}
which is a contradiction to $u_{0,2}(\xi+\eta)=u_{0,2}(\xi)+u_{0,2}(\eta)$.
\end{rem}
\begin{rem} If the assumption of relevancy is dropped, we must start with a time consistent system of utility functions $u_{0,2},u_{0,1},u_{1,2}$. In that case we only have that $\Qr\ll\Pr$ and the result of the previous remark only holds $\Qr$ \as  
\end{rem}
\begin{rem} There is no reason that $u_{0,2}$ is additive on $L^\infty(\Fc_2)$ as the following example shows.  We take $\Omega=[0,1]\times [0,1]$, $\Fc_2$ is the product sigma algebra of the Borel sigma algebas on $[0,1]$, the measure $\Pr$ is the product measure of the usual Lebesgue measures. $\Fc_0$ is the trivial sigma algebra and $\Fc_1$ is generated by the first coordinate mapping.  For $\xi\in L^\infty(\Fc_2),\xi\ge 0$ we define
$$
u_{0,2}(\xi)=\int_0^1d\alpha\int_0^\infty dx\,\Pr[\xi(\alpha,.)\ge x]^{1+\alpha}.
$$
\end{rem}
\section{A continuous time result}
In this section we use a filtration indexed by the time interval $[0,T]$.  This filtration $\(\Fc_t\)_{0\le t\le T}$ does not necessarily fulfil the usual assumptions.  The only assumption is that $\Fc_T$ is generated by $\cup_{0\le t< T}\Fc_t$.  We also suppose that a family of coherent utility functions $u_{t,s},0\le t\le s\le T$, $u_{t,s}\colon L^\infty(\Fc_s)\rightarrow L^\infty(\Fc_t)$ is given.  We assume the following time consistency:  for $t\le s\le v$ we have $u_{t,v}=u_{t,s}\circ u_{s,v}$. 
\begin{theo} We assume the notation introduced in this section.  We suppose that for $0\le t < T$, the sigma algebra $\Fc_T$ is atomless conditionally to $\Fc_t$. If $u_{0,T}$ is relevant, Lebesgue continuous and commonotone then there is a probability $\Qr\sim\Pr$ such that for all $\xi\in L^\infty(\Fc_T)$: $u_{0,T}(\xi)=\Er_\Qr[\xi]$.
\end{theo}
{\bf Proof }  The results of the previous section show that on each $L^\infty(\Fc_t)$, the utility function $u_{0,T}$ is linear.  The utility function $u_{0,T}$ is therefore linear on the vector space $\cup_{t<T}L^\infty(\Fc_t)$. This space is sequentially dense in $L^\infty(\Fc_T)$ for the 
Mackey topology (simply use the martingale convergence theorem). Because of Lebesgue continuity, the utility function $u_{0,T}$ is therefore linear on $L^\infty(\Fc_T)$.  It is therefore given by a probability measure $\Qr\ll\Pr$.  But since the utility function is relevant we find that $\Qr\sim\Pr$.
\begin{rem} The previous results apply for most filtrations used in finance and insurance.  For instance filtrations of Brownian Motion in one oir several dimensions, filtrations generated by most L\'evy processes and so on.  In other words {\it commonotonicity and time consistency are not good friends}.
\end{rem}

\end{document}